\documentclass[11pt,reqno,a4paper]{amsart}

\usepackage{array}
\usepackage{enumitem}
\usepackage{mathtools}
\newcommand\scalemath[2]{\scalebox{#1}{\mbox{\ensuremath{\displaystyle #2}}}}
\usepackage{pgfplots}
\usepackage{soul}

\usepackage{amsmath,amssymb,amsfonts}
\usepackage{tikz}
\usepackage{hyperref}
\usepackage{appendix}
\usepackage{scalerel}
\usepackage{graphicx}
\usepackage{subcaption}
\usepackage{amsthm}
\usepackage{multirow}
\usepackage{pdflscape}
\usepackage{afterpage}
\usepackage{capt-of} 
\usepackage{lipsum}
\usepackage{booktabs}
\usepackage{siunitx}
\usepackage{esvect}

\tikzset{decorated arrows/.style={
    postaction={
        decorate,
        decoration={
            markings,
            mark=between positions 0 and 1 step 15mm with {\arrow[black]{stealth};}
            }
        },
    }
}
 
\tikzset{decorated arrows2/.style={
    postaction={
        decorate,
        decoration={
            markings,
            mark=at position 15mm with {\arrow[black]{stealth};}
            }
        },
    }
}

\usepgfplotslibrary{colormaps}
\usetikzlibrary{arrows}

\makeatletter
\def\set@curr@file#1{%
  \begingroup
    \escapechar\m@ne
    \xdef\@curr@file{\expandafter\string\csname #1\endcsname}%
  \endgroup
}
\def\quote@name#1{"\quote@@name#1\@gobble""}
\def\quote@@name#1"{#1\quote@@name}
\def\unquote@name#1{\quote@@name#1\@gobble"}
\makeatother
\usepackage{graphics}

\newtheorem{theorem}{Theorem}
\newtheorem{lemma}{Lemma}
\newtheorem{corollary}{Corollary}

\newcommand{\RN}[1]{%
  \textup{\uppercase\expandafter{\romannumeral#1}}%
}
      
\author[J. P. S. M. de Carvalho]{Jo\~ao P.S. Maur\'icio de Carvalho$^{*1,2}$ \\
\\
$^*$\MakeLowercase{jocarvalho@fc.up.pt}  \\
\\
$^1$Faculty of Sciences, University of Porto, \\ Rua do Campo Alegre s/n, Porto 4169-007, Portugal \\
\\
$^2$Centre for Mathematics, University of Porto, \\ Rua do Campo Alegre s/n, Porto 4169-007, Portugal
}

\begin{document}

\thanks{JPSMC was supported by CMUP, Portugal (UIDP/MAT/00144/2020), which is funded by
Funda\c{c}\~ao para a Ci\^encia e a Tecnologia (FCT).\\ $^*$Corresponding author.}

\title[Dynamics of SARS-CoV-2 infection]
{A novel approach to understanding CoViD-19: exploring the interplay of SARS-CoV-2 and CTL response}

\date{\today}

\begin{abstract}

Facing a global challenge with over 6.9 million fatalities, severe acute respiratory syndrome coronavirus 2 (SARS-CoV-2), the causative agent of CoViD-19, demands novel and comprehensive approaches to understand its complex dynamics. This paper introduces a non-integer order model, capturing the intricate interplay between SARS-CoV-2 and the host's cytotoxic T lymphocytes (CTLs) response. Our work reveals a unique parameter space, in which an endemic state of SARS-CoV-2 and a CTL response-free equilibrium can coexist -- a crucial finding in our quest to decipher this pervasive virus. We further explore the basic reproduction number, assessing how different model parameters can potentially inhibit or fuel the infection's progression. Through extensive numerical simulations, we scrutinize the impact of varying the order of the fractional derivative and employing diverse CTL proliferation functions. This study significantly enriches our understanding of CoViD-19 immunopathology, offering invaluable insights that could guide future research and therapeutic strategies.
\end{abstract}

\maketitle

%
%
%
%
%
%
%

\section{Introduction}

\label{intro}
The first official case of CoViD-19 in China emerged in December 2019 and it was associated with the Huanan seafood market, in Wuhan \cite{Hu2021}. Since then SARS-CoV-2 has reached 231 countries and territories having infected more than 765 million people, thus becoming a global pandemic and causing more than 6.9 million deaths worldwide \cite{Vitiello2020,WHO2021,Worldometers2021}. The most frequently noticed symptoms in infected people are fever, cough and respiratory disorders \cite{Vitiello2020,Sun2020}. The most affected people are elderly and adults over 60 years old and/or those with comorbidities, such as obesity, diabetes, oncological diseases, heart problems, among others \cite{RothanByrareddy2020,Moriconi2020,Bajgain2021}.

Today, about three years after the first outbreak, there are already several vaccines that prevent coronavirus infection. Unfortunately vaccines are still not a fast enough method to fight the pandemic, especially in economically low-income countries \cite{Yi2020,Tracey2021}. Furthermore, the period of immunity that vaccines confer on people is not yet known exactly. As a consequence, an increase in the number of infected people has occurred in some countries with a large percentage of vaccinated individuals \cite{Taylor2021}.

It is becoming increasingly important to try to understand how the immune system reacts to SARS-CoV-2 in order to find alternatives while non-priority, young and economically disadvantaged people are not vaccinated \cite{Vitiello2020,Cunha2020,SetteCrotty2021}. Many mathematicians have proposed several models describing the dynamics of CoViD-19 in the population \cite{MauricioPinto2021,Elie2020,Kochanczyk2020,Liu2021,Cooper2020}. However, as far as mathematical modelling is concerned, existing studies on dynamical systems involving SARS-CoV-2 and the immune system have not yet been very thorough. Wang \emph{et al}.~\cite{Wang2020} developed a mathematical model to understand the dynamics and impact that SARS-CoV-2 has on target cells and on the body's immune response. The numerical results of this model allowed the authors to conclude that anti-inflammatory treatments are an asset in the recovery of infected individuals. In addition, this type of treatment helps to decrease the viral load, especially when it reaches its peak. The same results were obtained when the treatment involved antiviral drugs. Chatterjee \emph{et al}.~\cite{ChatterjeeBasir2020} designed a mathematical model in order to analyse how a certain drug administered at regular intervals helps the recovery of infected individuals. The authors concluded that if the drug is applied in the right dosage and adjusted to each individual's condition, it is effective in combating SARS-CoV-2. Also, Bairagi \emph{et al}.~\cite{BairagiAdak2017} proposed a model for the dynamics of the immune response to an infection. The authors formulated four functions to model the immune response:

\medskip

\begin{itemize}

\item[] $\boldsymbol{f_1(I,C) = qIC}$
\smallskip

\noindent {\it ``CTL proliferation depends both on infected cells density and CTL population"}; \\

\item[] $\boldsymbol{f_2(I,C) = qI}$
\smallskip

\noindent {\it ``CTL production is assumed to depend on infected cells density only"}; \\

\item[] $\boldsymbol{f_3(I,C) = \dfrac{qIC}{\varepsilon C + 1}}$
\smallskip

\noindent {\it ``CTL expansion saturates as the number of CTL 
grows to relatively high numbers. The level at which CTL expansion saturates is expressed in the parameter $\varepsilon$"}; \\

\item[] $\boldsymbol{f_4(I,C) = \dfrac{qI}{a + \varepsilon I}}$
\smallskip

\noindent {\it ``saturated type CTL production rate"}, where $a$ is 
the half-saturation constant,
\end{itemize}

\medskip

\noindent and $I(t)$ and $C(t)$ represent the population of infected cells and CTL, respectively. The proliferation rate of CTL is given by $q$. In this paper, we will analyse $f_1$, $f_2$, $f_3$ and $f_4$ adapted to SARS-CoV-2 infection dynamics and we will call them \emph{CTL proliferation functions}.

\subsection*{Advantages of fractional order models}

Fractional order (FO) models have been increasingly used and consulted in the literature in recent years due to their capacity in describing non-linear dynamics \cite{Muresan2016,Sweilam2017,Valerio2013}. These models have a great advantage over models of integer order equations since non-integer order models take into account a larger number of degrees of freedom. As such, the understanding about their dynamics and behaviours becomes more realistic. In recent times, there has been a growing interest to model epidemiological systems via non-integer order equations and in the last two and a half years, by virtue of circumstances, this interest has intensified in modelling CoViD-19 \cite{Shah2020,Awais2020,Oud2021}.


\subsection*{Model development and goals}

Motivated by all the reasons mentioned above, we built a FO model for the dynamics of SARS-CoV-2 in the presence of the immune response modelled by four CTL proliferation functions. This work has three main goals: (i) to examine the role of the order of fractional derivative $\alpha$, on the efficacy of the immune response, (ii) to explore the immune response for distinct CTL proliferation functions and $\alpha$ values, in the presence of SARS-CoV-2, (iii) find a space of parameters for which a SARS-CoV-2 endemic and a CTL response-free equilibria can coexist. In (ii), the numerical simulations of the model were obtained through the subroutine developed by Diethelm and Freed (1999) in \cite{DiethelmFreed1999}. A huge advantage of machine language is that it allows us to obtain numerical solutions very close to the real solutions, which are often extremely difficult or even impossible to obtain analytically. However, numerical methods always give us approximate solutions and not exact ones, since there is an error associated with each iteration.

\subsection*{Scope and adaptability of the model}

While the specific focus of this study is the SARS-CoV-2 virus, with model parameters chosen accordingly, it's important to note that the structure and methodologies of the model are not exclusive to this particular pathogen. Our model could feasibly be adapted to analyse other viral infections by adjusting the relevant parameters. Thus, this model can serve as a \emph{toy model}, providing a flexible framework for understanding the interactions between host immune responses and a variety of viral pathogens, not limited to SARS-CoV-2.

\subsection*{Structure}

In Section \ref{descricao} we describe the proposed model and analyse some of its properties. In Section \ref{EQUILIBRIA} we computed the equilibria of the model, studied the \emph{basic reproduction number} $\mathcal{R}_0$, and found a parameter space for which a SARS-CoV-2 endemic and a CTL response-free equilibria can coexist. In Section \ref{RESULTS} we outline some simulations of our dynamical system for different parameter values and for four proliferation functions. Finally, in Section \ref{conclusions} we draw the conclusions about our work.

\section{The model}
\label{descricao}
We analyse a FO model subdivided into four compartments: target/healthy cells susceptible to infection $T(t)$, infected cells $I(t)$, SARS-CoV-2 $V(t)$, and CTL $C(t)$ (see Table \ref{variables}). Let

$$
\Lambda = \left\{ (\lambda,\beta,\mu,k,\delta,N,c,q,\sigma, \varepsilon, a) \in (\mathbb{R}^+)^{11} \right\}
$$

\medskip

\noindent be the set of parameters of our model. The proliferation rate of healthy cells is given by $\lambda^{\alpha}$. Healthy cells are infected by the virus at a rate $\beta^{\alpha}$. The natural death rate of healthy and infected cells are given by $\mu^{\alpha}$ and $\delta^{\alpha}$, respectively. The term $N \delta^{\alpha} I$ represents the production of new viruses by infected cells during their lifetime. The viruses are cleared from the body at a rate of $c^{\alpha}$. CTL are produced through the proliferation functions $f_n(I, C)$, where $f_1 = q^{\alpha}IC$, $f_2 = q^{\alpha}I$, $f_3 = q^{\alpha}IC/(\varepsilon C + 1)$ and  $f_4 = q^{\alpha}I/(a + \varepsilon I)$ \cite{BairagiAdak2017}. We consider $q^{\alpha}$ to be the proliferation rate of CTL. The level of saturation of CTL expansion is given by $\varepsilon^{\alpha}$. Parameter $a^{\alpha}$ is the half-saturation constant. The natural death rate of CTL is given by $\sigma^{\alpha}$. Also, $0 < \alpha \leq 1$ is the order of the fractional derivative. The description and value of these parameters can be found in Table \ref{tabela}. The system of FO equations is

\begin{equation}\label{modelo}
\begin{array}{lcl}
\dfrac{d^{\alpha}T}{dt^{\alpha}} = \lambda^{\alpha} - \beta^{\alpha}VT - \mu^{\alpha} T \\
\\
\dfrac{d^{\alpha}I}{dt^{\alpha}} = \beta^{\alpha}VT - k^{\alpha}IC - \delta^{\alpha}I \\
\\
\dfrac{d^{\alpha}V}{dt^{\alpha}} = N \delta^{\alpha}I -c^{\alpha}V \\
\\
\dfrac{d^{\alpha}C}{dt^{\alpha}} = f_n (I,C) - \sigma^{\alpha} C .
\end{array}
\end{equation}

\medskip

We apply the principle of a FO derivative proposed by Caputo, \emph{i.e.}

\begin{equation*}
\label{caputo}
\dfrac{d^{\alpha}y(t)}{{dt}^{\alpha}}  =  I^{p-\alpha} y^{(p)} (t), \,\,\,\,\, t>0,
\end{equation*}

\noindent where $p=[\alpha]$ is the $\alpha$-value rounded up to the nearest integer, $y^{(p)}$ is the $p$\,-th derivative of $y(r)$ and $I^{p_1}$ is the Riemann-Liouville operator (Please see \cite{Samko1993} and references therein)

\begin{equation*}
\label{liouville}
I^{p_1} z(t)  =  \dfrac{1}{\Gamma(p_1)} \displaystyle \int_0^t (t - t')^{p_1 - 1} z(t') dt' ,
\end{equation*}

\noindent where $\Gamma(\cdot)$ is the gamma function.

\begin{table}[ht!]
\centering
\def\arraystretch{1.3} 
\begin{tabular}{lc}
\hline\noalign{\smallskip}
\textbf{Variable} & \textbf{Symbol} \\
\noalign{\smallskip}\hline\noalign{\smallskip}
Target cells & $T(t)$   \\
Infected cells & $I(t)$  \\
SARS-CoV-2 & $V(t)$  \\
CTL & $C(t)$  \\
\noalign{\smallskip}\hline
\end{tabular}
\caption{\label{variables}Description of the variables of model (\ref{modelo}).}
\end{table}


\subsection{Model properties analysis}
The solutions of the system (\ref{modelo}) remain non-negative for the entire domain, $t>0$. Let $\mathbb{R}_+^4 = \left\{x \in \mathbb{R}^4 \,\, | \,\, x\geq0 \right\}$ and $x(t) = \left[T(t), I(t), V(t), C(t) \right]^T$. First, we quote the following Generalized Mean Value Theorem \cite{OdibatShawagfeh2007} and corollary.

\begin{lemma}
\cite{OdibatShawagfeh2007} Assume $f(x) \in C[a,b]$ and $D_a^{\alpha}f(x) \in C[a,b]$, where $0<\alpha\leq1$. Thus

\begin{equation*}
\label{lema1}
	f(x)  =  f(a)+\dfrac{1}{\Gamma(\alpha)}(D_a^{\alpha}f)(\varepsilon) \cdot (x-a)^{\alpha}
\end{equation*}

\medskip

\noindent for $a \leq \varepsilon \leq x, \forall x \in (a,b]$. 
\end{lemma}

\begin{corollary} \label{corolario}
Suppose that $f(x) \in C[a,b]$ and $D_a^{\alpha}f(x) \in C(a,b]$, for $0<\alpha \leq 1$.

\medskip

\begin{enumerate}
\item If $D_a^{\alpha}f(x) \geq 0$, $\forall x \in (a,b)$, then $f(x)$ is non-decreasing for each $x \in [a,b]$; \label{remains}
\\
\item If $D_a^{\alpha}f(x) \leq 0$, $\forall x \in (a,b)$, then $f(x)$ is non-increasing for each $x \in [a,b]$.
\end{enumerate}
\end{corollary}

\noindent This proves the main theorem.

\begin{theorem}
There is a solution $x(t) = \left[T(t), I(t), V(t), C(t) \right]^T$ to the system (\ref{modelo}) in the entire domain $(t \geq 0)$ and it is unique. Furthermore, the solution remains in $\mathbb{R}_+^4$.
\end{theorem}

\begin{proof}
In \cite[Theorem 3.1, Remark 3.2]{Lin2007}, we can observe that the solution of the initial value problem exists and is unique in $\mathbb{R}_0^+$. To this end, it is enough to prove that the non-negative orthant $\mathbb{R}_+^4$ is positively invariant. So, we have to demonstrate that the vector field points to $\mathbb{R}_+^4$ in each hyperplane, thus limiting the non-negative orthant. Hence, for model (\ref{modelo}), we get:

\begin{equation*}
\label{limitada}
\begin{array}{lcl}
\dfrac{d^{\alpha}T}{dt^{\alpha}} \, \big|_{T=0}   =  \lambda^{\alpha} \geq 0 \\
\\
\dfrac{d^{\alpha}I}{dt^{\alpha}} \, \big|_{I=0}  =  \beta^{\alpha}VT \geq 0  \\
\\
\dfrac{d^{\alpha}V}{dt^{\alpha}} \, \big|_{V=0}  =  N \delta^{\alpha}I  \geq 0 \\
\\
\dfrac{d^{\alpha}C}{dt^{\alpha}} \, \big|_{C=0}  =  f_n (I,C)  \geq 0 .
\end{array}
\end{equation*}

\medskip

\noindent We concluded, by \emph{(\ref{remains})} of Corollary \ref{corolario}, that the solution will remain in $\mathbb{R}_+^4$.
\end{proof}

\section{Equilibria and basic reproduction number}
\label{EQUILIBRIA}
Throughout this section, we will be using $f_n (I,C) = f_1(I,C)$. We will study the following equilibria of the model (\ref{modelo}): (i) disease-free equilibrium point; (ii) CTL response-free equilibrium point; (iii) SARS-CoV-2 endemic equilibrium point and compute the \emph{basic reproduction number}. Furthermore, we will find a space of parameters for which a SARS-CoV-2 endemic and a CTL response-free equilibria can coexist.

\subsection{Disease-free equilibria}
A disease-free equilibrium point of the model (\ref{modelo}) is obtained via imposing $I = V = C = 0$. Let 

\begin{equation*}
\label{DFeq_point}
X_0 = (T_0, I_0, V_0, C_0) = \left(\dfrac{\lambda^{\alpha}}{\mu^{\alpha}}, 0, 0, 0 \right) .
\end{equation*}

\medskip

\noindent be the disease-free equilibrium point where there is no infection present in the organism. The linearization matrix of (\ref{modelo}) at a general point $X = (T, I, V, C) \in (\mathbb{R}_0^+)^4$ is given by:
 
\medskip 
 
\begin{equation}
\label{jacob_matrix}
\begin{array}{lcl}
J(X)=\left(\begin{array}{cccc}
-\left(\beta^{\alpha} V + \mu^{\alpha}\right) & 0 & -\beta^{\alpha} T & 0 \\ 
\\
\beta^{\alpha} V & -\left(k^{\alpha} C + \delta^{\alpha}\right) & \beta^{\alpha} T & -k^{\alpha} I \\
\\
0 & N \delta^{\alpha} & - c^{\alpha} & 0 \\
\\
0 & q^{\alpha} C & 0 & q^{\alpha}I - \sigma^{\alpha}
\end{array}\right)
\end{array}.
\end{equation}

\subsection{Basic reproduction number}
In a cellular scenario, the \emph{basic reproduction number} $\mathcal{R}_0$, is the number of secondary infections due to a single infected cell in a susceptible healthy cell population \cite{DriesscheWatmough2002}. Through \cite[Lemma 1]{DriesscheWatmough2002} the \emph{basic reproduction number} can be obtained by the following expression:

\begin{equation}
\label{R0_0}
\mathcal{R}_0 = \rho \left(FV^{-1}\right) ,
\end{equation}

\medskip

\noindent where $F$ is the matrix of the new infection terms, $V$ is the matrix of the remaining terms and $\rho$ is the spectral radius of the matrix $FV^{-1}$. For the model \eqref{modelo} we get:

\medskip

\begin{equation*}
\label{FV2}
\begin{array}{lcl}
F=\left(\begin{array}{cc}
0 & \beta^{\alpha} T_0 \\ 
\\
0 & 0
\end{array}\right) , \quad
V=\left(\begin{array}{cc}
\delta^{\alpha} & 0 \\ 
\\
-N \delta^{\alpha} & c^{\alpha}
\end{array}\right) .
\end{array}
\end{equation*}

\medskip

\noindent Then, using \eqref{R0_0} we compute $\mathcal{R}_0$ as:

\medskip

\begin{equation}
\label{R0}
\mathcal{R}_0 = \dfrac{\beta^{\alpha} T_0 N \delta^{\alpha}}{(k^{\alpha} C_0 + \delta^{\alpha} )c^{\alpha}} =   \dfrac{\beta^{\alpha} \lambda^{\alpha}}{\mu^{\alpha} c^{\alpha}} N .
\end{equation}

\bigskip

\subsection{Stability of the disease-free equilibrium $\boldsymbol{X_0}$}

We are now able to understand the asymptotic behaviour of $X_0$ by relating it to the \emph{basic reproduction number} $\mathcal{R}_0$. By Theorem 2 of \cite{DriesscheWatmough2002} we obtain Lemma \ref{lemma_X0}.

\begin{lemma} \label{lemma_X0}
If $\mathcal{R}_0 < 1$, then the disease-free equilibrium point $X_0$ is locally asymptotically stable. If $\mathcal{R}_0 > 1$, then $X_0$ is unstable. 
\end{lemma}

\begin{proof}
We know that at the disease-free equilibrium point $X_0$, the matrix \eqref{jacob_matrix} takes the following form:

\medskip

\begin{equation}
\label{jacob_matrix_X0}
\begin{array}{lcl}
J(E)=\left(\begin{array}{cccc}
-\mu^{\alpha} & 0 & -\dfrac{\beta^{\alpha} \lambda^{\alpha}}{\mu^{\alpha}} & 0 \\ 
\\
0 & -\delta^{\alpha} & \dfrac{\beta^{\alpha} \lambda^{\alpha}}{\mu^{\alpha}} & 0 \\
\\
0 & N \delta^{\alpha} & - c^{\alpha} & 0 \\
\\
0 & 0 & 0 & - \sigma^{\alpha}
\end{array}\right)
\end{array},
\end{equation}

\medskip

\noindent and its eigenvalues are given by:

\medskip

$\gamma_1 = -\sigma^{\alpha},$ \quad $\gamma_2 = - \mu^{\alpha},$

\bigskip

$
\gamma_3 = \dfrac{-\mu^{\alpha}\left(c^{\alpha}+\delta^{\alpha}\right) + \displaystyle\sqrt{4N\beta^{\alpha}\delta^{\alpha}\lambda^{\alpha}\mu^{\alpha}+\left(c^{\alpha} \mu^{\alpha}\right)^2 - 2c^{\alpha}{\mu^2}^{\alpha} \delta^{\alpha} + \left(\delta^{\alpha} \mu^{\alpha}\right)^2}}{2\mu^{\alpha}},
$

\bigskip

\noindent and

\bigskip

$
\gamma_4 = \dfrac{-\mu^{\alpha}\left(c^{\alpha}+\delta^{\alpha}\right) - \displaystyle\sqrt{4N\beta^{\alpha}\delta^{\alpha}\lambda^{\alpha}\mu^{\alpha}+\left(c^{\alpha} \mu^{\alpha}\right)^2 - 2c^{\alpha}{\mu^2}^{\alpha} \delta^{\alpha} + \left(\delta^{\alpha} \mu^{\alpha}\right)^2}}{2\mu^{\alpha}}.
$

\bigskip

\noindent It is easy to verify that $\gamma_1<0$, $\gamma_2<0$ and  $\gamma_4<0$. If

\begin{equation*}\label{eigenvalues_negative}
\begin{array}{lcl}
&& \gamma_3 < 0 \\
\\
&\Leftrightarrow& -\mu^{\alpha}\left(c^{\alpha}+\delta^{\alpha}\right) + \displaystyle\sqrt{4N\beta^{\alpha}\delta^{\alpha}\lambda^{\alpha}\mu^{\alpha}+\left(c^{\alpha} \mu^{\alpha}\right)^2 - 2c^{\alpha}{\mu^2}^{\alpha} \delta^{\alpha} + \left(\delta^{\alpha} \mu^{\alpha}\right)^2} < 0 \\
\\
&\Leftrightarrow& 4N\beta^{\alpha}\delta^{\alpha}\lambda^{\alpha}\mu^{\alpha}+\left(c^{\alpha} \mu^{\alpha}\right)^2 - 2c^{\alpha}{\mu^2}^{\alpha} \delta^{\alpha} + \left(\delta^{\alpha} \mu^{\alpha}\right)^2 < \left(c^{\alpha} \mu^{\alpha}\right)^2 + 2c^{\alpha}{\mu^2}^{\alpha} \delta^{\alpha} + \left(\delta^{\alpha} \mu^{\alpha}\right)^2 \\
\\
&\Leftrightarrow& 4N\beta^{\alpha}\delta^{\alpha}\lambda^{\alpha}\mu^{\alpha} - 4 c^{\alpha}{\mu^2}^{\alpha} \delta^{\alpha} < 0 \\
\\
&\Leftrightarrow& \dfrac{\beta^{\alpha} \lambda^{\alpha} }{ \mu^{\alpha}c^{\alpha} }N < 1 \\
\\
&\overset{\eqref{R0}}{\Leftrightarrow}& \mathcal{R}_0 < 1,
\end{array}
\end{equation*}

\bigskip

then all eigenvalues have negative real part and $X_0$ is locally asympotically stable.
\end{proof}

\subsection{CTL response-free and SARS-CoV-2 endemic equilibria}

The following results show us that there is a region where a CTL response-free and SARS-CoV-2 endemic equilibria can coexist.

\begin{lemma} \label{LEMMA1}
If \,$\mathcal{R}_0 > 1$, then the CTL response-free equilibrium point exists.
\end{lemma}

\begin{proof}
If we impose $C=0$ in the model (\ref{modelo}) we get the following equilibrium point:

\begin{equation*}
\label{CTLeq_point}
\begin{array}{lcl}
X_1 &=& (T_1,I_1,V_1,C_1) \\
\\
&=& \left(\dfrac{c^{\alpha}}{N \beta^{\alpha}}, \dfrac{N \beta^{\alpha} \lambda^{\alpha} - c^{\alpha} \mu^{\alpha}}{N \beta^{\alpha} \delta^{\alpha}},  \dfrac{N \beta^{\alpha} \lambda^{\alpha} - c^{\alpha} \mu^{\alpha}}{\beta^{\alpha} c^{\alpha}}, 0\right) .
\end{array}
\end{equation*}

\medskip

\noindent It is clear that all cells and virus populations are non-negative. In this case it is easy to verify that $T_1 > 0$ and $C_1 = 0$. In order to $I_1$ and $V_1$ we have $(N \beta^{\alpha} \lambda^{\alpha} - c^{\alpha} \mu^{\alpha})/(N \beta^{\alpha} \delta^{\alpha}) > 0$ and $(N \beta^{\alpha} \lambda^{\alpha} - c^{\alpha} \mu^{\alpha})/(\beta^{\alpha} c^{\alpha}) > 0$, respectively. So, if

\begin{equation}
\label{CTLeq_point2}
N \beta^{\alpha} \lambda^{\alpha} - c^{\alpha} \mu^{\alpha}  >  0  \quad \Leftrightarrow \quad N \beta^{\alpha} \lambda^{\alpha}  >  c^{\alpha} \mu^{\alpha} \quad \overset{\eqref{R0}}{\Leftrightarrow} \quad \mathcal{R}_0  >  1 , 
\end{equation}

\medskip

\noindent then $I_1, V_1 > 0$. Therefore, the CTL response-free equilibrium point exists.
\end{proof}

Now, we set the following real constant will be useful in the upcoming results:

\begin{equation}
\label{varphi0}
\varphi_0 = 1 + A,
\end{equation}

\medskip

\noindent where $A = \sigma^{\alpha} N \beta^{\alpha} \delta^{\alpha}/(q^{\alpha} c^{\alpha} \mu^{\alpha})$. Since $A > 0$, then $\varphi_0 > 1$.

\begin{lemma}
If \,$\mathcal{R}_0 > \varphi_0$, then the SARS-CoV-2 endemic and CTL response-free equilibria can coexist.
\end{lemma}

\begin{proof}
Similar to Lemma \ref{LEMMA1}, we compute the endemic equilibrium point of system (\ref{modelo}) and we get

\begin{equation*}
\label{ENDeq_point}
X_2 = (T_2, I_2, V_2, C_2),
\end{equation*}

\noindent where

\begin{equation*}
\label{ENDeq_point2}
\begin{array}{lcl}
T_2 &=& \dfrac{\lambda^{\alpha} c^{\alpha} q^{\alpha}}{\sigma^{\alpha} N \beta^{\alpha} \delta^{\alpha} + c^{\alpha} \mu^{\alpha} q^{\alpha}} \\
\\
I_2 &=&  \dfrac{\sigma^{\alpha}}{q^{\alpha}} \\
\\
V_2 &=& \dfrac{N \delta^{\alpha} \sigma^{\alpha}}{q^{\alpha} c^{\alpha}} \,\,\, = \,\,\, \dfrac{N \delta^{\alpha}}{c^{\alpha}}I_2 \\
\\
C_2 &=& B \Big[ q^{\alpha}\left(N \beta^{\alpha} \lambda^{\alpha} - c^{\alpha} \mu^{\alpha}\right) - \sigma^{\alpha} N \beta^{\alpha} \delta^{\alpha} \Big] ,
\end{array}
\end{equation*}

\medskip

\noindent and $B = \delta^{\alpha}/\left[ (\sigma^{\alpha} N \beta^{\alpha} \delta^{\alpha} + c^{\alpha} \mu^{\alpha} q^{\alpha} ) k^{\alpha}\right] > 0$. It is clear that $T_2,I_2,V_2 > 0$. Since $B > 0$, if

\begin{eqnarray}
\nonumber && q^{\alpha}\left(N \beta^{\alpha} \lambda^{\alpha} - c^{\alpha} \mu^{\alpha}\right) - \sigma^{\alpha} N \beta^{\alpha} \delta^{\alpha} > 0 \\
\nonumber \\
\nonumber \Leftrightarrow && q^{\alpha} N \beta^{\alpha} \lambda^{\alpha} - q^{\alpha} c^{\alpha} \mu^{\alpha} - \sigma^{\alpha} N \beta^{\alpha} \delta^{\alpha} > 0 \\
\nonumber \\
\nonumber \Leftrightarrow && q^{\alpha} \left( N \beta^{\alpha} \lambda^{\alpha} - c^{\alpha} \mu^{\alpha} \right) > \sigma^{\alpha} N \beta^{\alpha} \delta^{\alpha}  \\
\nonumber \\
\nonumber \Leftrightarrow && N \beta^{\alpha} \lambda^{\alpha} - c^{\alpha} \mu^{\alpha} > \dfrac{\sigma^{\alpha} N \beta^{\alpha} \delta^{\alpha}}{q^{\alpha}} \\
\nonumber \\
\nonumber \Leftrightarrow && \dfrac{N \beta^{\alpha} \lambda^{\alpha} - c^{\alpha} \mu^{\alpha}}{c^{\alpha} \mu^{\alpha}} > \dfrac{\sigma^{\alpha} N \beta^{\alpha} \delta^{\alpha}}{q^{\alpha} c^{\alpha} \mu^{\alpha}} \\
\nonumber \\
\overset{\eqref{R0}}{\Leftrightarrow} && \mathcal{R}_0 > 1 + \dfrac{\sigma^{\alpha} N \beta^{\alpha} \delta^{\alpha}}{q^{\alpha} c^{\alpha} \mu^{\alpha}} \label{R01mais} \\
\nonumber \\
\overset{\eqref{varphi0}}{\Leftrightarrow} && \mathcal{R}_0 > \varphi_0 \label{R0phi0} ,
\end{eqnarray}

\medskip

\noindent then $C_2 > 0$. Hence, condition \eqref{R0phi0} provides the existence of the SARS-CoV-2 endemic equilibrium point. 
Moreover, looking at expressions (\ref{CTLeq_point2}) and (\ref{R01mais}), we can easily verify that the CTL response-free and SARS-CoV-2 endemic equilibria can coexist. Figure \ref{coexistencia} illustrates this coexistence.
\end{proof}


\begin{center}
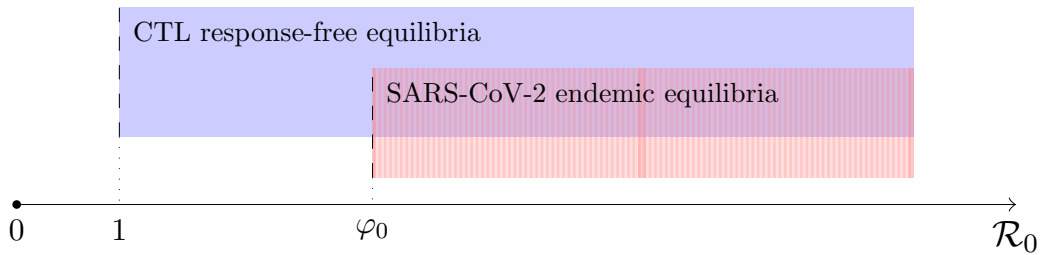
\begin{figure}[ht!]
\scalemath{1.8}{
\begin{tikzpicture}

\path[left color=blue!20, right color=blue!20] (1.25,0.50) -- (1.25,1.45) -- (7.05,1.45) -- (7.05,0.50) -- cycle;
\path[left color=red!50, right color=red!50, opacity=0.25] (3.10,0.20) -- (3.10,1.0) -- (7.05,1.0) -- (7.05, 0.20) -- cycle;

\draw [very thin, dashed] (1.25,0.50)--(1.25,1.45) node[below right=0.9, scale = 0.53] {\textcolor{black}{CTL response-free equilibria}}; 
\draw [very thin, dashed] (3.10,0.2)--(3.10,1.0) node[below right=0.9, scale = 0.53]{\textcolor{black}{SARS-CoV-2 endemic equilibria}}; 

\draw [very thin, dotted] (1.25,0.5)--(1.25,0) node[below=0.9, scale = 0.6]{\textcolor{black}{$1$}}; 
\draw [very thin, dotted] (3.10,0.2)--(3.10,0) node[below=0.9, scale = 0.6]{\textcolor{black}{$\varphi_0$}}; 

\filldraw[black] (0.5,0) circle (0.7pt) node[below=0.075em, scale = 0.6] {0}; 

\draw [->, very thin] (0.5,0)--(7.8,0) node[below=0.075em, scale = 0.7]{$\mathcal{R}_0$}; 

\end{tikzpicture}
}
\caption{Coexistence space of SARS-CoV-2 endemic and CTL response-free equilibria for $\mathcal{R}_0 > \varphi_0$.}
\label{coexistencia} 
\end{figure}
\end{center}

\section{Numerical simulations and results}
\label{RESULTS}
In this section we simulate the variation of $\mathcal{R}_0$ as a function of $\beta^{\alpha}$, $\mu^{\alpha}$, $c^{\alpha}$, $N$ and $\lambda^{\alpha}$ (see \eqref{R0}) and sketched the dynamics of model (\ref{modelo}) for different values of $\alpha \in [0,1]$ and CTL proliferation functions. The parameter values used in the numerical simulations are in Table \ref{tabela}. The initial conditions that we assume are $T(0) = 1000$, $I(0) = 0$, $V(0) = 10$ and $C(0) = 333$ \cite{CMAP2019}. Since the model of the authors of \cite{CMAP2019} is similar to ours, we consider that these initial conditions also fit our model. The numerical solutions were obtained by applying the subroutine of Diethelm and Freed (FracPECE subroutine) \cite{DiethelmFreed1999}.

\subsection*{Stability and convergence of the method}

The classical Adams-Bashforth-Moulton method \cite{ABM} used for first order differential equations is a good method and extremely useful in the sense that its stability properties provide safe information about in slightly stiff equations without glaring oscillations from rounding errors. We can find in \cite{DiethelmFreed1999} that these stability properties remain unchanged for differential equations of fractional order, that is, the stability of this method does not depend of the order of the fractional derivative $\alpha$. Of course, stability alone is not sufficient in practice to make sure that the numerical solution is a good approximation to the exact solution. As such, convergence is a problem that must also be analysed. The method we use, is a method that uses the same techniques as \cite{Linz1985} whose standard error is given by the expression

\bigskip

\begin{center}
$\max_j |y(t_j) - y_j| = O(h^2),$
\end{center}

\bigskip

\noindent where $h = \max_j (t_{j+1} - t_j)$ and $t_j \in [ t_0, t_0 + t^{\star} ]$ for some fixed $t^{\star} > 0$. The FracPECE method was already used by other authors in previous studies, such as \cite{MauricioPinto2021, JoCarvalhoPinto2022}. 

\begin{table}[ht!]
\centering
\scalebox{0.89}{
\def\arraystretch{1.6}
\begin{tabular}{ lclc }
\hline\noalign{\smallskip}
\textbf{Parameter} & \textbf{Symbol} & \textbf{Value} & \textbf{Reference}  \\
\noalign{\smallskip}\hline\noalign{\smallskip}
Proliferation rate of healthy cells & $\lambda^{\alpha}$ & $10$ cells mm$^{-3}$  & \cite{BairagiAdak2017}   \\
SARS-CoV-2 infection rate & $\beta^{\alpha}$ & $0.001$ virions mm$^3$ day$^{-\alpha}$  & \cite{Tang2017} \\
Bursting size for virus growth & $N$ & $10 - 2500$ virions cell$^{-1}$ & \cite{BairagiAdak2017} \\
Proliferation rate of CTL & $q^{\alpha}$ & $0.2$ & \cite{BairagiAdak2017} \\
Elimination rate of infected cells by CTL & $k^{\alpha}$ & $0.7$ cells & \cite{BairagiAdak2017} \\
Natural death rate of healthy cells & $\mu^{\alpha}$ & $0.01$ day$^{-\alpha}$ & \cite{BairagiAdak2017}  \\
Natural death rate of infected cells & $\delta^{\alpha}$ & $1$ day$^{-\alpha}$& \cite{Goncalves2020}  \\
Natural death rate of CTL & $\sigma^{\alpha}$ & $0.08$ cells day$^{-\alpha}$ & \cite{BairagiAdak2017} \\
Death rate of virus & $c^{\alpha}$ & $3$ day$^{-\alpha}$  & \cite{BairagiAdak2017} \\
Level of saturation of CTL expansion & $\varepsilon$ & $0.01$ cells  & \cite{BairagiAdak2017} \\
Half-saturation constant & $a$ & $120$ cells  & \cite{BairagiAdak2017} \\
\noalign{\smallskip}\hline
\end{tabular}
}
\caption{\label{tabela}Parameter values used in numerical simulations of model (\ref{modelo}).}
\end{table}

\begin{figure*}
\includegraphics[width=1.022\textwidth]{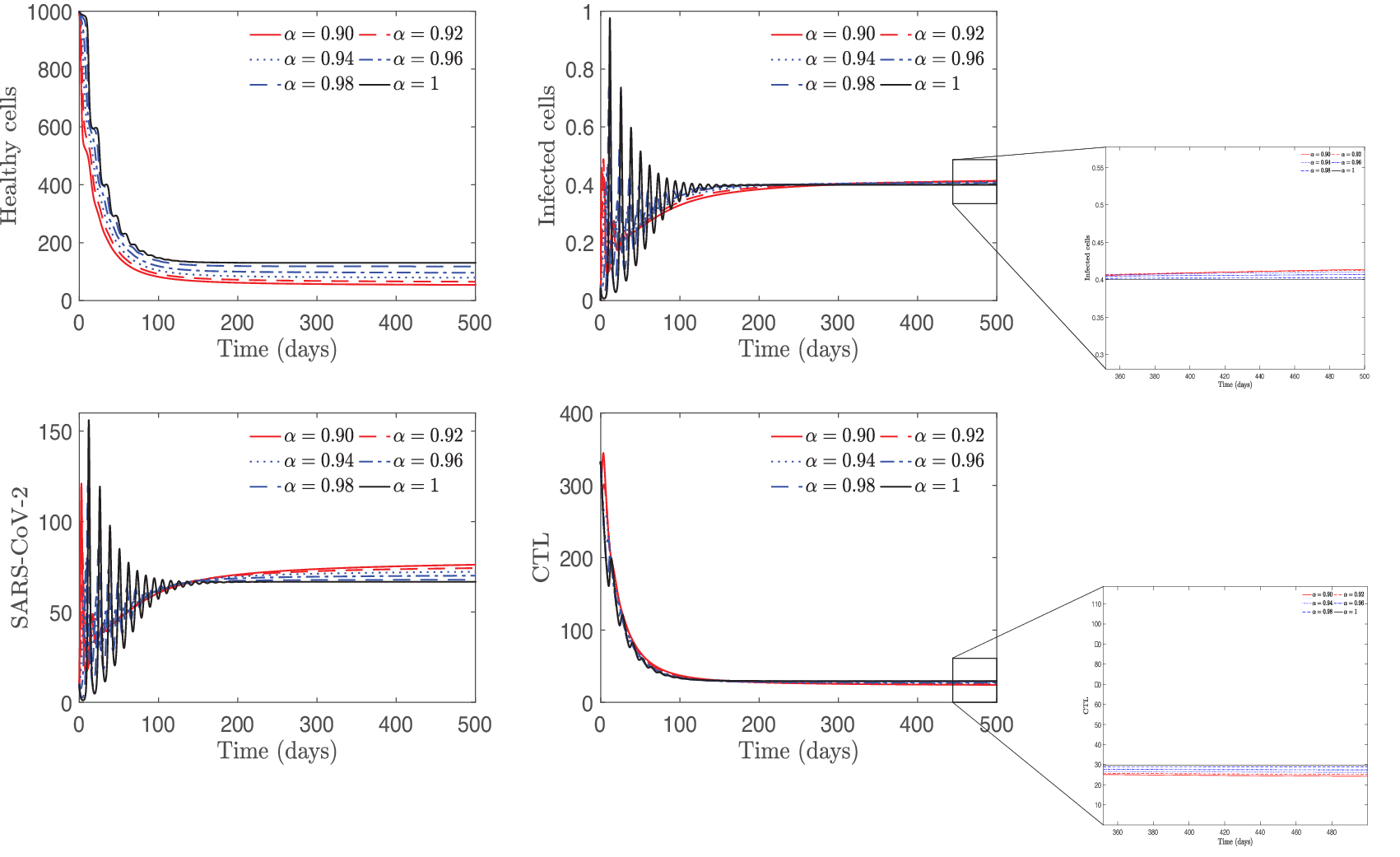}
\caption{Dynamic behaviour of $T$, $I$, $V$ and $C$ populations for $\alpha \in \{0.90, 0.92, 0.94, 0.96, 0.98, 1 \}$. We consider $f_n (I,C) = f_1(I,C)$. Initial conditions and parameter values are given in the text and in the Table \ref{tabela}, respectively.}
\label{pop_behaviour}    
\end{figure*}

Figure \ref{pop_behaviour} shows the behaviour of the four classes for different values of the fractional order derivative $\alpha$. It is concluded that for lower $\alpha$ values there is a lower number of healthy T cells and CTL. However, and regarding the CTL, in the first 150 days, for low values of $\alpha$ there is a higher number of CTL. Also in the first 150 days, higher values of $\alpha$ relate to a lower number of CTL. We can also conclude that the lower the $\alpha$ value, the higher the viral load of SARS-CoV-2 in the organism. Regarding to infected cells, it is observed that the higher the $\alpha$ value, the lower the number of these cells. 

\begin{figure*}
\includegraphics[width=0.48\textwidth]{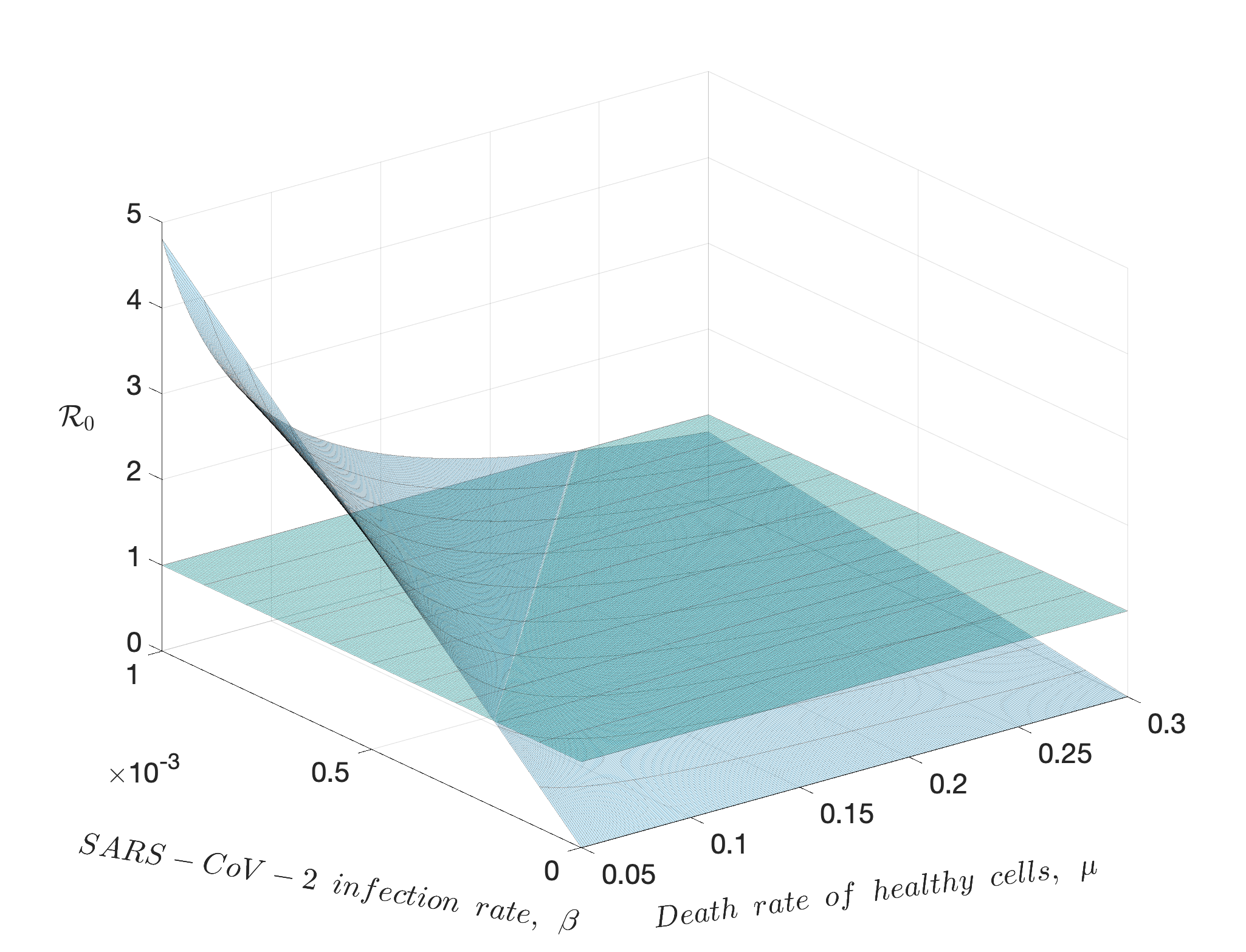}
\includegraphics[width=0.48\textwidth]{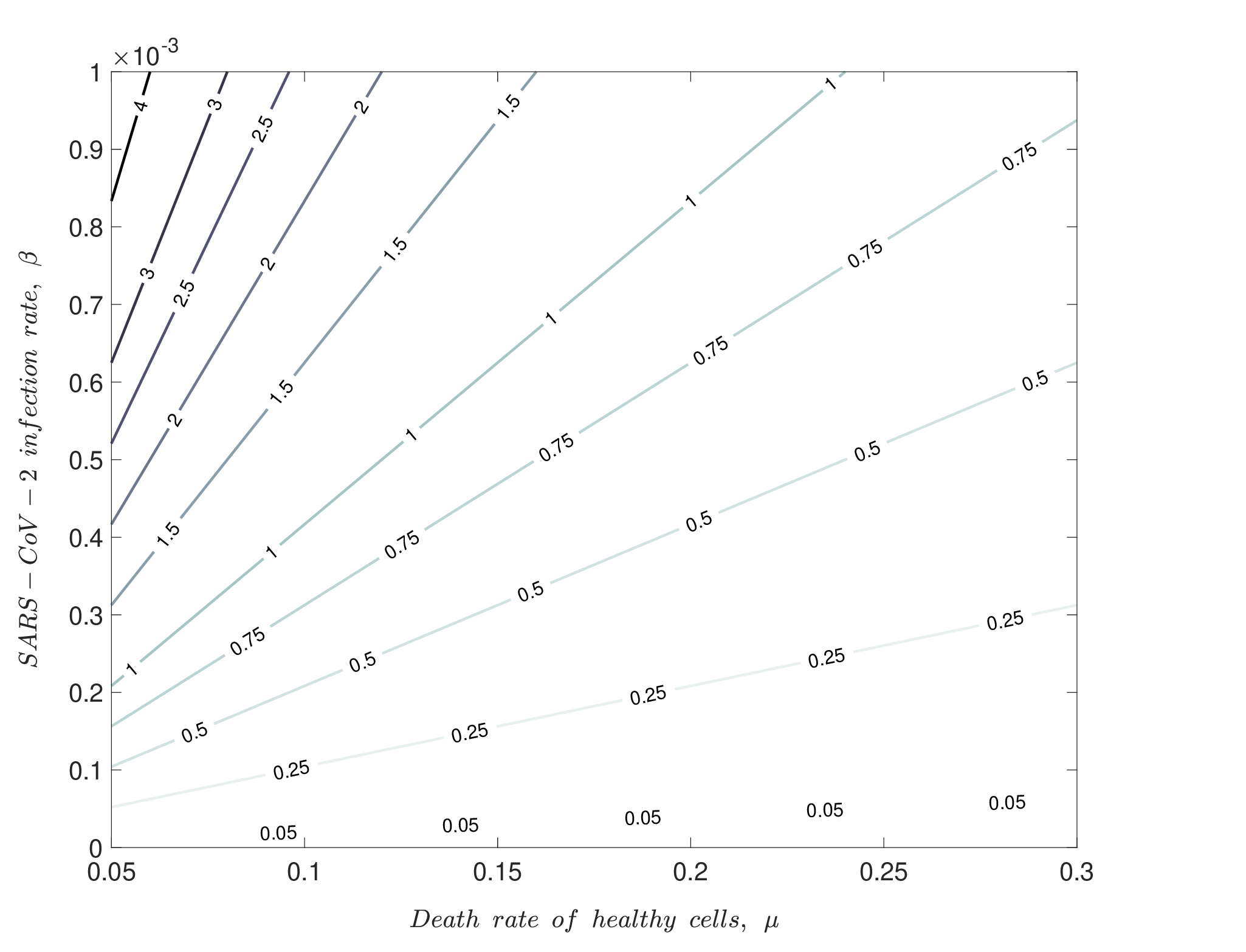}
\caption{Effect of SARS-CoV-2 infection rate $\beta$, and death rate of healthy cells $\mu$, on the \emph{basic reproduction number} $\mathcal{R}_0$. Initial conditions and parameter values are given in the text and in the Table \ref{tabela}, respectively, except for $\beta$ and $\mu$. We consider $\alpha = 1$.}
\label{R0_variation1}    
\end{figure*}

\begin{figure*}
\includegraphics[width=0.48\textwidth]{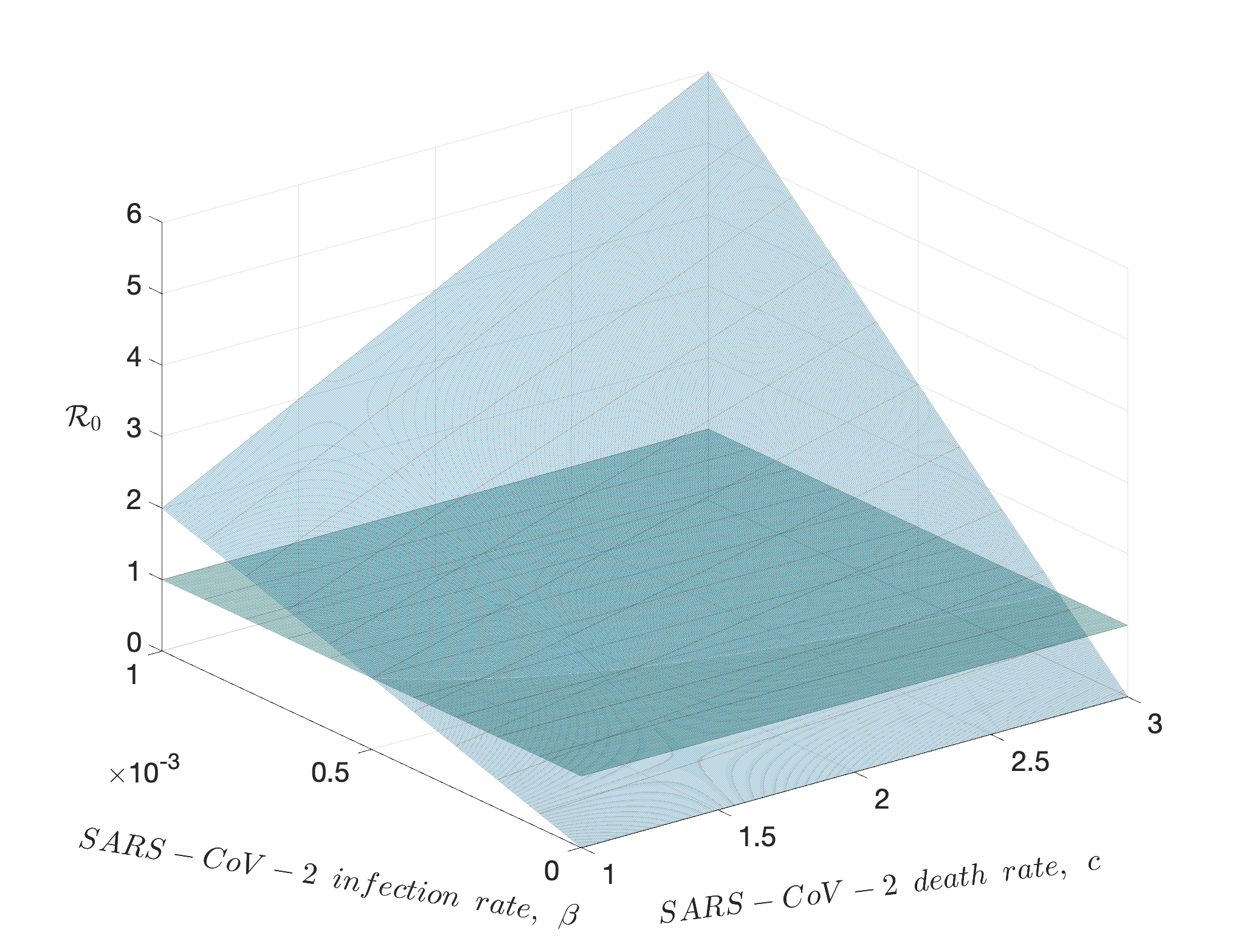}
\includegraphics[width=0.48\textwidth]{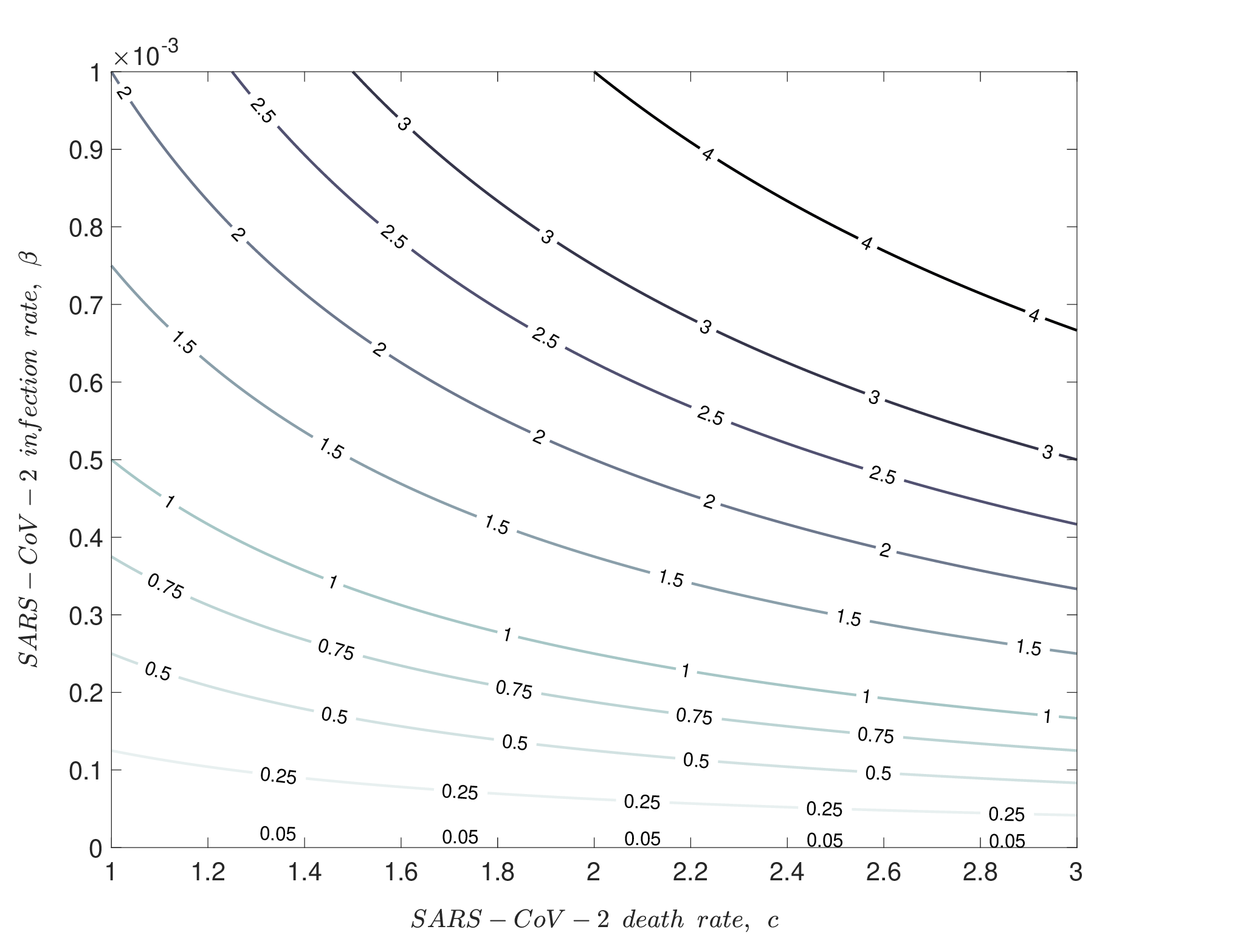}
\caption{Effect of SARS-CoV-2 infection rate $\beta$ and SARS-CoV-2 death rate $c$, on the \emph{basic reproduction number} $\mathcal{R}_0$. Initial conditions and parameter values are given in the text and in the Table \ref{tabela}, respectively, except for $\beta$ and $c$. We consider $\alpha = 1$.}
\label{R0_variation2}  
\end{figure*}

\begin{figure*}
\includegraphics[width=0.48\textwidth]{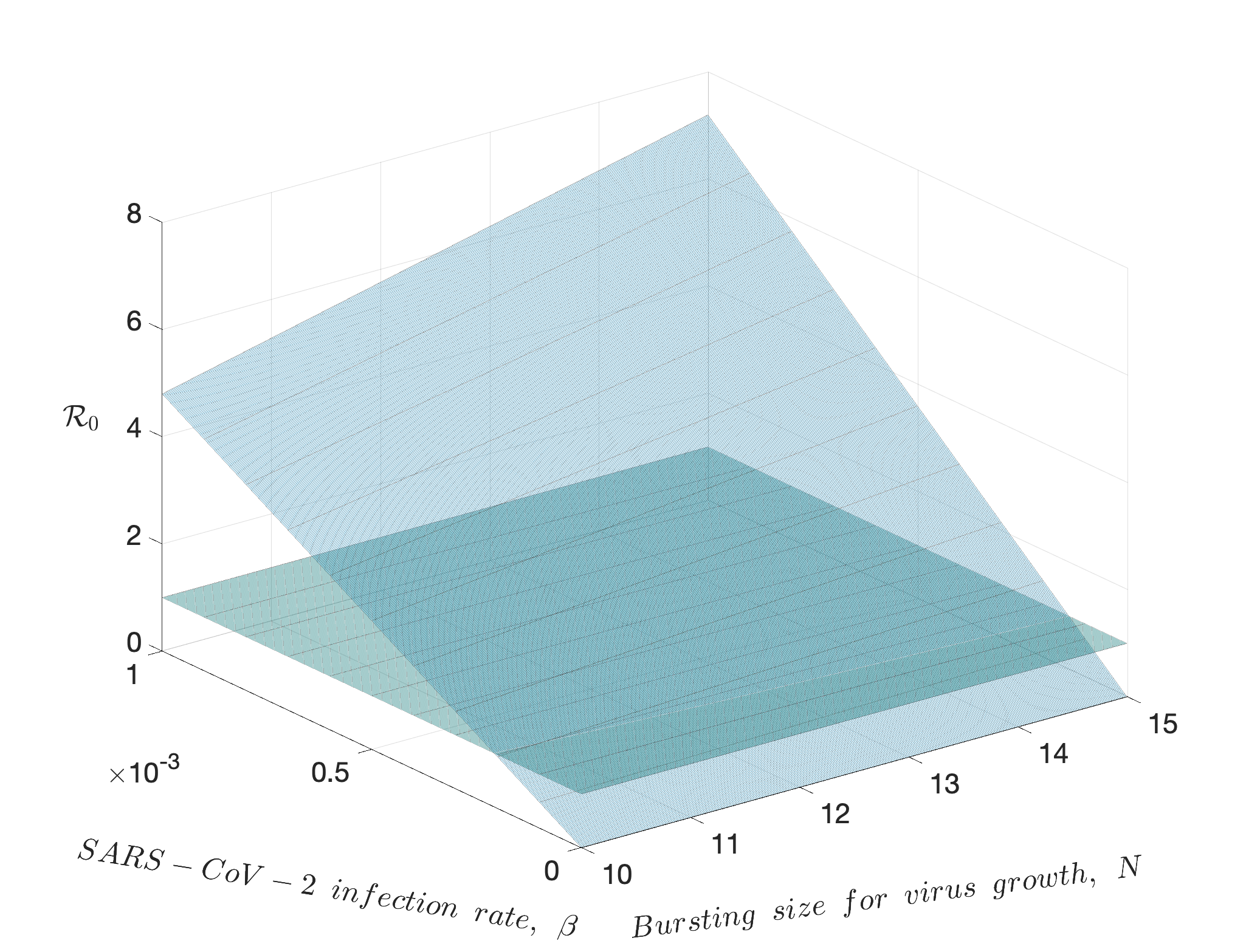}
\includegraphics[width=0.48\textwidth]{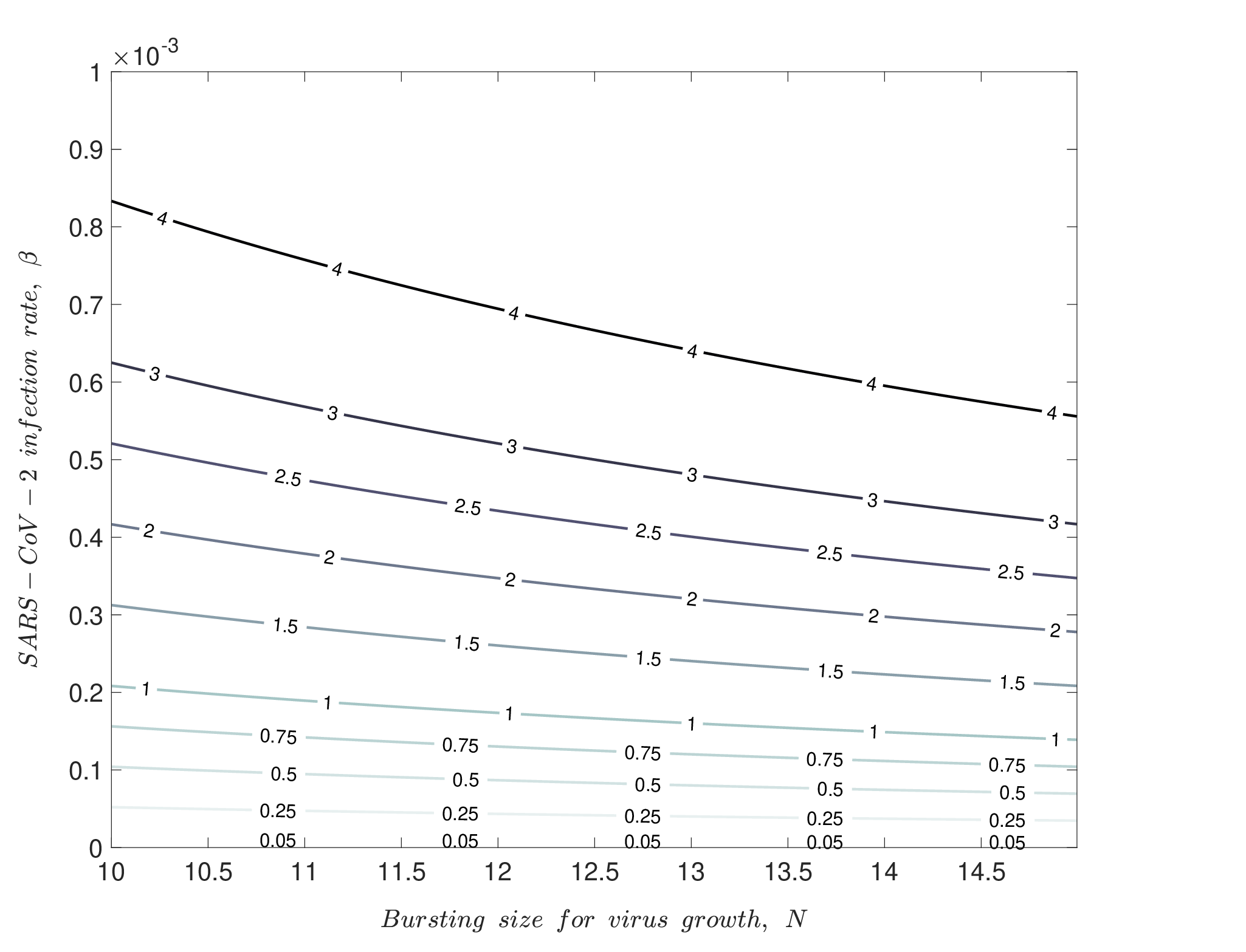}
\caption{Effect of SARS-CoV-2 infection rate $\beta$, and bursting size for virus growth $N$, on the \emph{basic reproduction number} $\mathcal{R}_0$. Initial conditions and parameter values are given in the text and in the Table \ref{tabela}, respectively, except for $\beta$ and $N$. We consider $\alpha = 1$.}
\label{R0_variation3}
\end{figure*}

\begin{figure*}
\includegraphics[width=0.48\textwidth]{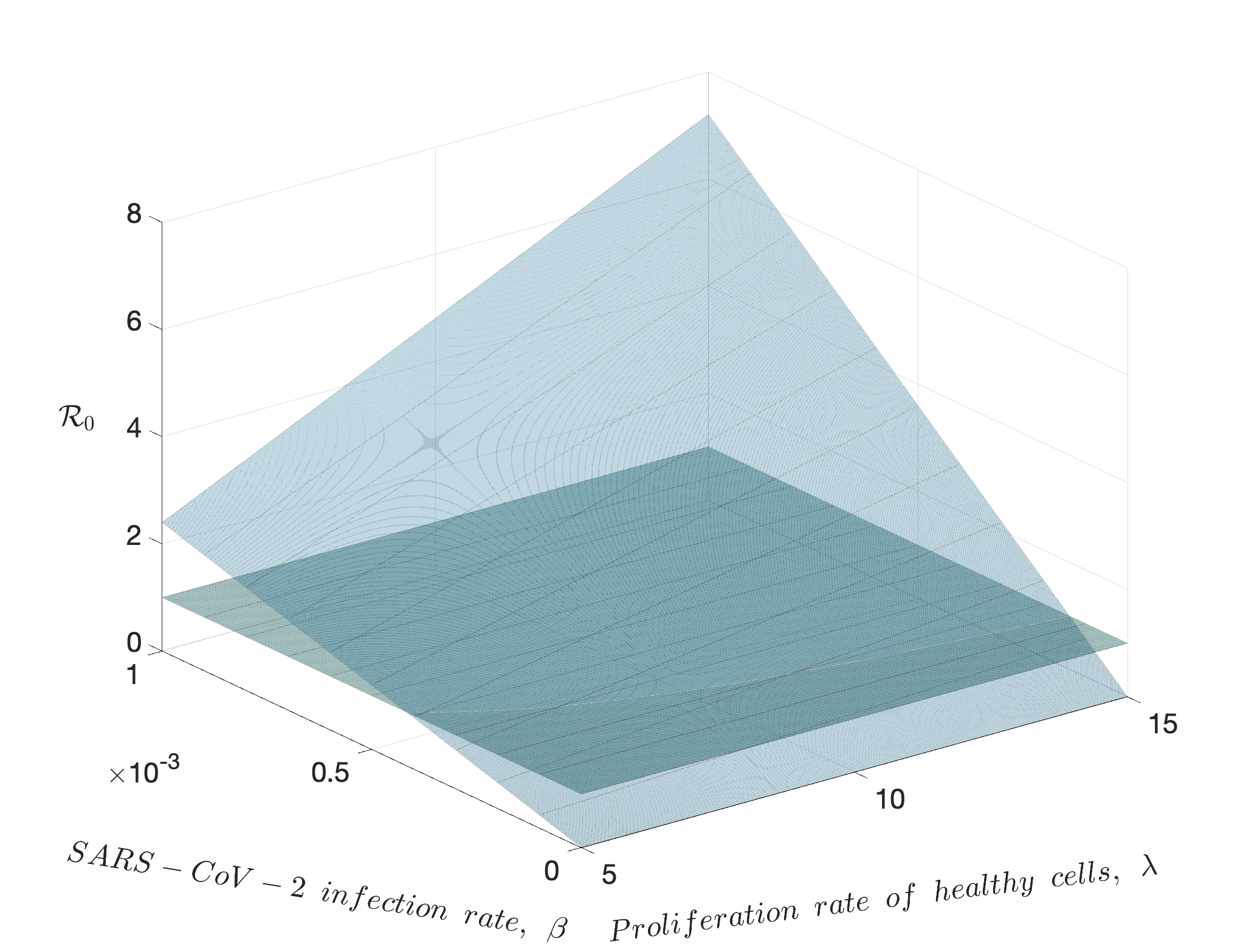}
\includegraphics[width=0.48\textwidth]{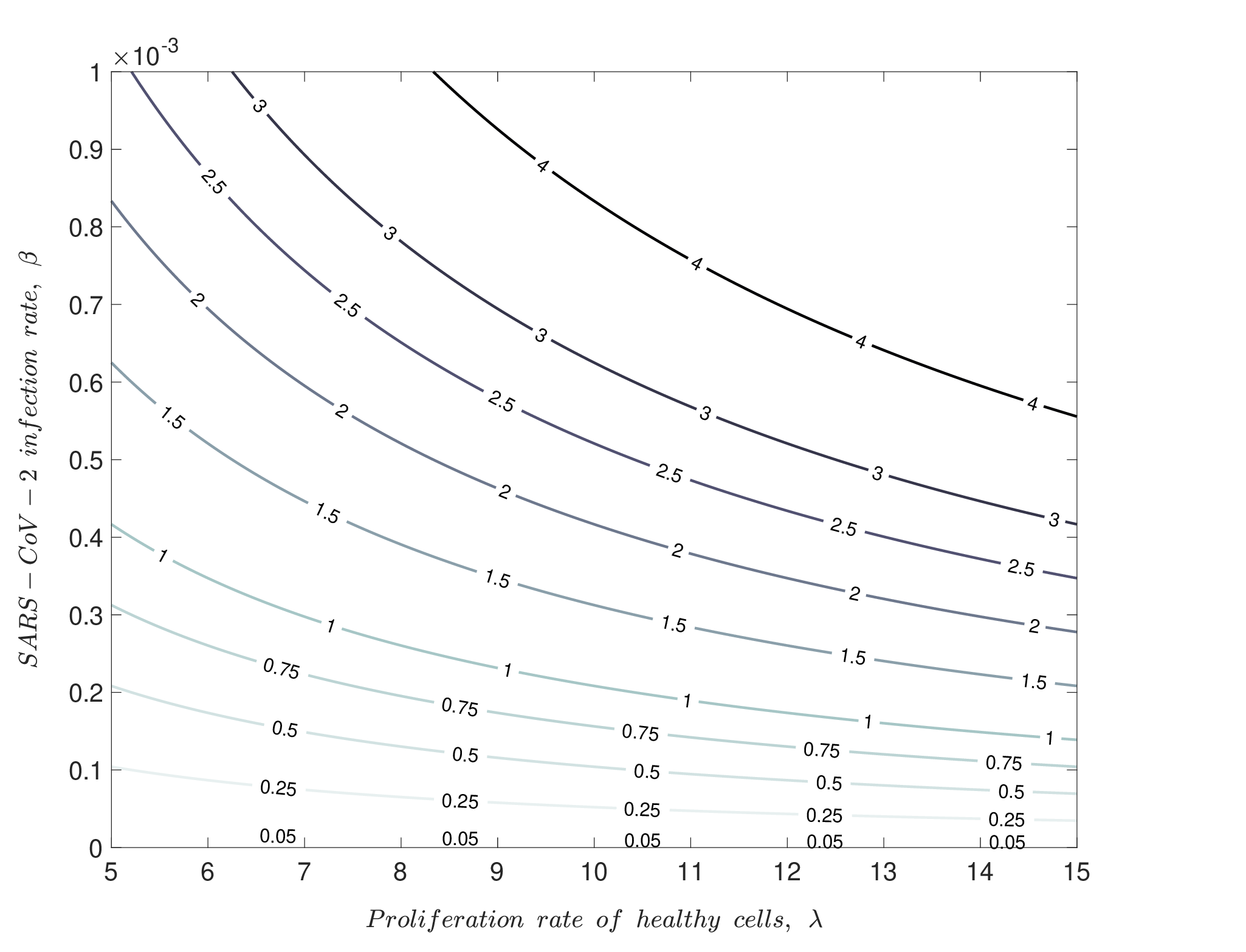}
\caption{Effect of SARS-CoV-2 infection rate $\beta$, and proliferation rate of healthy cells $\lambda$, on the \emph{basic reproduction number} $\mathcal{R}_0$. Initial conditions and parameter values are given in the text and in the Table \ref{tabela}, respectively, except for $\beta$ and $\lambda$. We consider $\alpha = 1$.}
\label{R0_variation4}
\end{figure*}

Figures \ref{R0_variation1}~--~\ref{R0_variation4} show the effect of SARS-CoV-2 infection rate $\beta^{\alpha}$, death rate of healthy cells $\mu^{\alpha}$, SARS-CoV-2 death rate $c^{\alpha}$, bursting size for virus growth $N$, and proliferation rate of healthy cells $\lambda^{\alpha}$, on the value of $\mathcal{R}_0$, for $\alpha = 1$.

From literature in general, when $\mathcal{R}_0 < 1$ then the disease tends to slow down until it disappears. Otherwise, when $\mathcal{R}_0 > 1$, the disease spreads further and further through the population \cite{DriesscheWatmough2002}. In all figures we can observe that $\mathcal{R}_0$ increases with the value of $\beta$, promoting the progression of SARS-CoV-2 infection. We can interpret the variation in this parameter value with the data from the World Health Organization (WHO) \cite{WHO2021}. For example, during the periods when the governments of the countries allowed the population to be in lockdown, naturally the number of contacts increased and consequently the number of cases also increased. With this, at a cellular level, the viral load is higher and the transmission rate of SARS-CoV-2 between cells also increase. This behaviour was reflected in an increase of $\mathcal{R}_0$, as we can see in Figures \ref{R0_variation1}~--~\ref{R0_variation4}.

However, in Figure \ref{R0_variation1}, for $\beta < 0.2$ or for $\mu > 0.25$, the \emph{basic reproduction number} is less than 1, which leads us to conclude that the infection will eventually disappear. Figure \ref{R0_variation2} suggests that low values of $\beta$ combined with high values of $c$ result in a small value of $\mathcal{R}_0$. We observed by Figure \ref{R0_variation3} that $N$ has no apparent significant influence of the progression of SARS-CoV-2 infection. Furthermore, for $\beta < 0.15$ the infection tends to disappear. In Figure \ref{R0_variation4} we can see that $\mathcal{R}_0$ has lower values for combinations of low values of $\beta$ and high values of $\lambda$.

\begin{figure}
\includegraphics[width=1.14\textwidth]{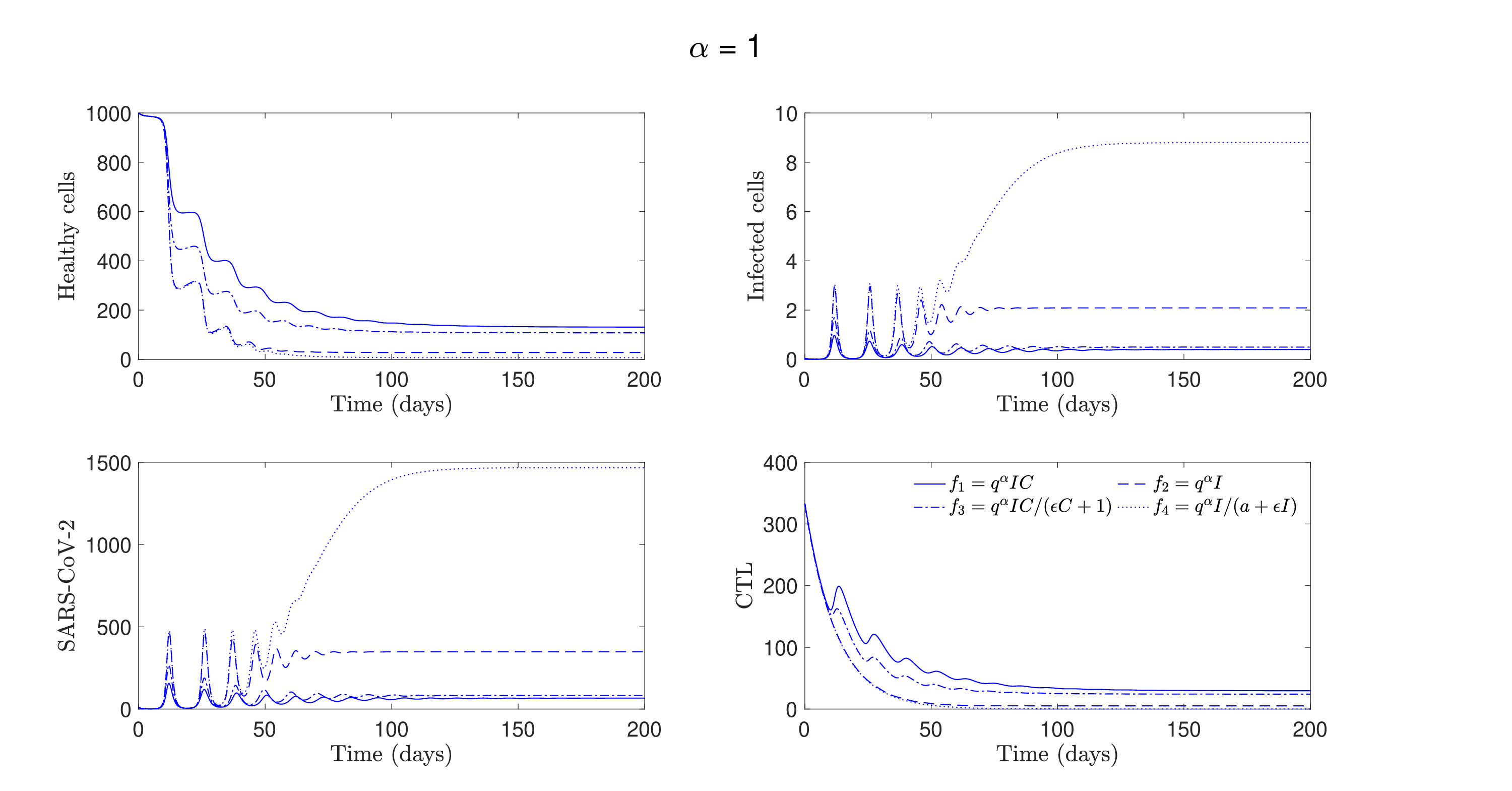}
\caption{Dynamics of model (\ref{modelo}) for $f_1$, $f_2$, $f_3$ and $f_4$. Initial conditions and parameter values are given in the text and in the Table \ref{tabela}, respectively. We consider $\alpha = 1$.}
\label{comportamento1}       
\end{figure}

\begin{figure}
\includegraphics[width=1.14\textwidth]{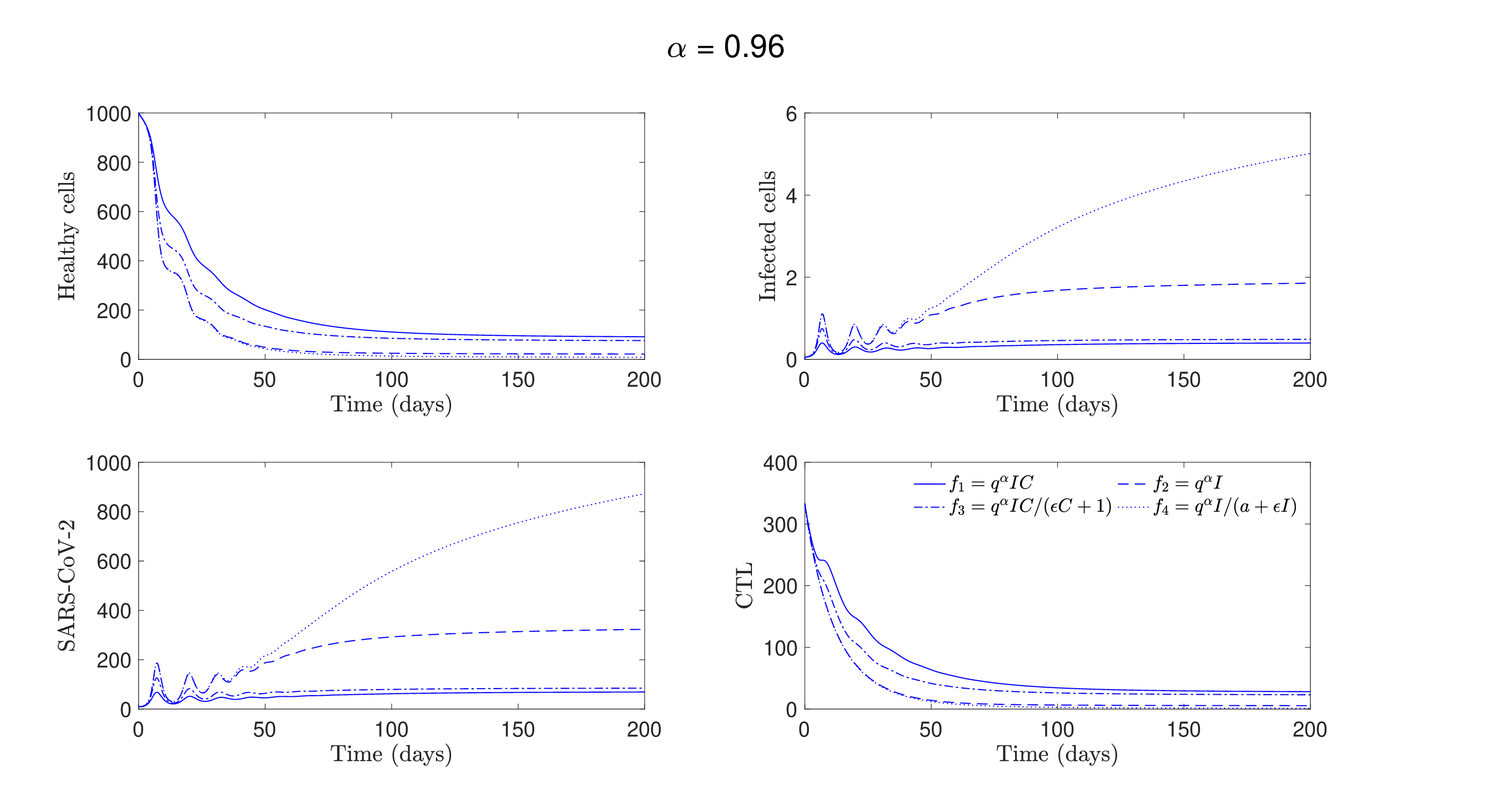}
\caption{Dynamics of model (\ref{modelo}) for $f_1$, $f_2$, $f_3$ and $f_4$. Initial conditions and parameter values are given in the text and in the Table \ref{tabela}, respectively. We consider $\alpha = 0.96$.}
\label{comportamento99}
\end{figure}

\begin{figure}
\includegraphics[width=1.14\textwidth]{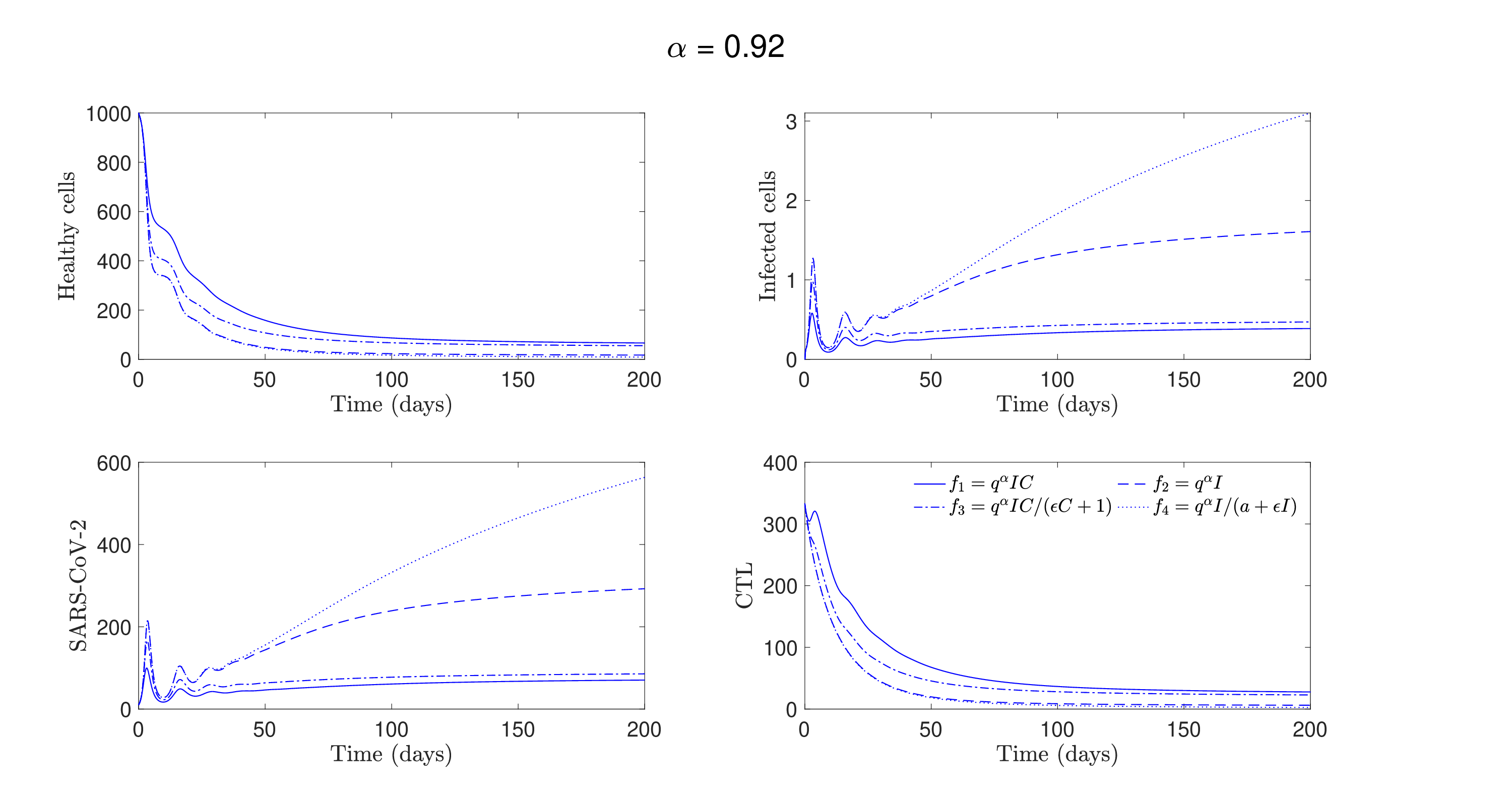}
\caption{Dynamics of model (\ref{modelo}) for $f_1$, $f_2$, $f_3$ and $f_4$. Initial conditions and parameter values are given in the text and in the Table \ref{tabela}, respectively. We consider $\alpha = 0.92$.}
\label{comportamento98}
\end{figure}

In Figures \ref{comportamento1}~--~\ref{comportamento98} we consider four different immune functions $f_n(I,C)$, specifically $f_1 (I,C) = q^{\alpha} I C$, $f_2 (I,C) =  q^{\alpha} I$, $f_3 (I,C) =  q^{\alpha} IC/(\varepsilon C +1)$ and $f_4 (I,C) =  q^{\alpha} I/(a + \varepsilon I)$, and different values of $\alpha$ (Please see figure legends). For all proliferation functions the system presents the same asymptotic behaviour converging to an endemic state. However, it can be observed that the number of infected cells and the SARS-CoV-2 viral load are higher when we consider the proliferation function $f_4$. It is also possible to notice that for the proliferation function $f_4$, there is a higher number of infected cells. This may be a consequence of the weak immune response (Please see CTL subplot). This weak immune response is also reflected in the viral load, which is also higher for $f_4$. The opposite is true for $f_1$. As the immune response is stronger, the number of infected cells and the viral load is lower. This dynamic occurs for all values of the order of the fractional derivative $\alpha$, with the particularity that the lower its value, the less severe is the epidemic state.

\section{Conclusions} \label{conclusions}

In this study, a FO model for the dynamics of SARS-CoV-2, responsible for CoViD-19, was analysed together with the human immune response, in particular CTL. The research essentially consisted of three key points:

\begin{itemize}
\item[{\bf (I)}] Theoretical analysis of the model to find a region where two equilibria can coexist;
\item[{\bf (II)}] Influence that the parameters of the model have on the progression of the virus within human organism;
\item[{\bf (III)}] Analysis of the immune response that human body is able to provide in the presence of SARS-CoV-2, for different $\alpha$ values.
\end{itemize}

Regarding the coexistence of two equilibria, we found that under particular conditions, namely for $\mathcal{R}_0 > \varphi_0$, the CTL response-free and SARS-CoV-2 endemic equilibria can coexist (see Figure \ref{coexistencia}). This result is quite important in that we find a range of values for the \emph{basic reproduction number} ($1 < \mathcal{R}_0 < \varphi_0$) for which there is no disease or immune response (which makes sense). However, for $\mathcal{R}_0 > \varphi_0$, the organism still has no immune response to the virus that is already present in the human body. This phenomenon should serve as a warning to health professionals and epidemiologists. 

One of the conclusions that could possibly have more impact in real life is related to Figure \ref{pop_behaviour}. We can see that for lower values of $\alpha$, the asymptotic behaviour of the populations stabilizes faster than for higher values of $\alpha$. This scenario can play an important role in predicting cell behaviour earlier than if we run simulations with just the integer order derivative ($\alpha = 1$).

With respect to the \emph{basic reproduction number} $\mathcal{R}_0$, simulations have revealed that when $\beta < 0.2$, then $\mathcal{R}_0 < 1$, which leads us to conclude that the infection tends to disappear. So human body will be healthy again and free of infection.

With regard to the analysis of the immune response to infection, we simulated the model for four CTL proliferation functions:

$$
f_1(I,C) = q^{\alpha}IC, \quad
f_2(I,C) = q^{\alpha}I, \quad
f_3(I,C) = \dfrac{q^{\alpha}IC}{\varepsilon C + 1}, \quad
f_4(I,C) = \dfrac{q^{\alpha}I}{a + \varepsilon I},
$$

\noindent using the subroutine of Diethelm and Freed \cite{DiethelmFreed1999}.

We concluded that the \emph{saturated type CTL production rate} $f_4$, is the function that most worsens the endemic state of infection, although all CTL proliferation functions asymptotically converge towards an endemic equilibrium. Furthermore, we verified that CTL proliferation functions that trigger stronger immune responses, such as $f_1$ and $f_3$, lead to fewer infected cells and a less pronounced viral load. All these conclusions about the immune response to SARS-CoV-2 are consistent for different values of $\alpha$. However, the lower the value of $\alpha$, the lower the endemic state of the infection. 

The results of the simulations are in line with what was expected. The method we used to find numerical solutions to our model is a method that is used in many published scientific papers, some of them on epidemic models similar to ours \cite{CMAP2019}. This leads us to believe that the accuracy of the numerical solutions of our model is quite reasonable. We address the readers who are interested in consulting the algorithm of Diethelm and Freed (1999) to reference \cite{DiethelmFreed1999}.

These models allow us to draw a lot of conclusions about how we should act in situations related to epidemics in different countries around the world. For example, these mathematical analysis could help in decision making by government leaders in combating possible economic losses that come from pandemic disasters, particularly due to SARS-CoV-2 \cite{MauricioPinto2021}.

Considering the scarce information relating the dynamics of SARS-CoV-2 and the immune system, we were able to obtain solid results that may contribute to the understanding of the disease and its mechanisms of action from a mathematical point of view.  We hope in the near future to present an improved version of this model, based on more biochemical information regarding parameter values. Additionally, we plan to include the effect of vaccination at a cellular level to study how the immune system is reinforced and what the impact is on other cells.

%
%
%

\bigskip

\section*{Declarations}

\subsection*{Ethical statement} The authors have no ethical statement to declare
\subsection*{Funding} JPSMC was supported by CMUP, Portugal (UIDP/MAT/00144/2020), which is funded by
Funda\c{c}\~ao para a Ci\^encia e a Tecnologia (FCT)
\subsection*{Conflict of Interest/Competing interests} The authors declare that they have no conflict of interest
\subsection*{Availability of data and material} This manuscript has no associated data
\subsection*{Code availability} Not applicable
\subsection*{Authors' contributions} JPSMC is the only author of this work

\end{document}